\documentclass[a4paper,11pt]{elsarticle}
\usepackage{amsfonts}
\usepackage{amstext}
\usepackage{amsthm}
\usepackage{amsmath}
\usepackage{paralist}
\usepackage{epsfig}
\usepackage{amssymb}
\usepackage{latexsym}
\usepackage{graphicx}
\usepackage[latin1]{inputenc}
\newcommand{\R}{\Bbb{R}}

\newtheorem{teor}{Theorem}[section]
\newtheorem{propo}{Proposition}[section]

\newtheorem{cor}{Corollary}[section]
\newtheorem{rem}{Remark}[section]

\newcommand{\n}{\noindent}

\newcommand {\fim}{\rule{0.5em}{0.5em}}

\begin{document}

\title{Positive solutions for non-variational fractional elliptic systems with negative exponents
}

\author{\textbf{Anderson L. A. de Araujo \footnote{\textit{E-mail addresses}:
anderson.araujo@ufv.br(A.
Araujo)}\footnote{A.L.A de Araujo was partially supported by FAPEMIG/FORTIS.}}\\ {\small\it Departamento de Matem\'{a}tica,
Universidade Federal de Vi\c{c}osa,}\\ {\small\it CCE, 36570-000, Vi\c{c}osa, MG, Brazil}\\
\textbf{Luiz F. O. Faria \footnote{\textit{E-mail addresses}:
luiz.faria@ufjf.edu.br(L. F. O. Faria)}\footnote{L. F.O Faria was partially supported by FAPEMIG CEX APQ 02374/17.}}\\ {\small\it Departamento de Matem\'{a}tica,
Universidade Federal de Juiz de Fora,}\\ {\small\it ICE, 36036-330, Juiz de Fora, MG, Brazil}\\
\textbf{Edir Junior F. Leite \footnote{\textit{E-mail addresses}:
edirjrleite@ufv.br (E.J.F. Leite)}}\\ {\small\it Departamento de Matem\'{a}tica,
Universidade Federal de Vi\c{c}osa,}\\ {\small\it CCE, 36570-000, Vi\c{c}osa, MG, Brazil}\\
\textbf{Ol\'{i}mpio H. Miyagaki \footnote{\textit{E-mail addresses}:
olimpio.hiroshi@ufjf.edu.br (O. H. Miyagaki)} \footnote{O. H. M. was supported in part by INCTMAT/BRAZIL and CNPQ/BRAZIL PROC. 307061/2018-3} }\\ {\small\it Departamento de Matem\'{a}tica,
Universidade Federal de Juiz de Fora,}\\ {\small\it ICE, 36036-330, Juiz de Fora, MG, Brazil}}

\date{}{

\maketitle

\markboth{abstract}{abstract}
\addcontentsline{toc}{chapter}{abstract}

\hrule \vspace{0,2cm}

\n {\bf Abstract}

In this paper, we study strongly coupled elliptic systems in non-variational form with negative exponents involving fractional Laplace operators. We investigate the existence, nonexistence, and uniqueness of the positive classical solution. The results obtained here are a natural extension of the results obtained by Ghergu, in \cite{ghergu}, for the fractional case.

\vspace{0.5cm}
\hrule\vspace{0.2cm}

\section{Introduction and main results}

The present paper deals with existence, nonexistence, and uniqueness of positive solutions for elliptic systems of the form

\begin{equation} \label{1}
\left\{
\begin{array}{llll}
(-\Delta)^{s}u = u^{-p}v^{-q} & {\rm in} \ \ \Omega\\
(-\Delta)^{t}v = u^{-r}v^{-\theta} & {\rm in} \ \ \Omega\\
u= v=0 & {\rm in} \ \ \R^n\setminus\Omega
\end{array}
\right.
\end{equation}
where $\Omega$ is a smooth bounded open subset of $\R^{n}$, $n\geq 2$, $0 <s,t < 1$, $r,q > 0$, $p, \theta\geq 0$ and the fractional Laplace operator $(-\Delta)^{s}$ is defined as

\[
(-\Delta)^{s}u(x) = C(n,s)\lim_{\varepsilon\searrow 0}\int\limits_{\R^{n}\setminus B_\varepsilon(x)}\frac{u(x)-u(y)}{\vert x-y\vert^{n+2s}}\; dy\, ,
\]
\n for all $x \in \R^{n}$ and

\[
C(n,s) = \left(\int\limits_{\R^{n}}\frac{1-\cos(\zeta_{1})}{\vert\zeta\vert^{n+2s}}\; d\zeta\right)^{-1}
\]

\n with $\zeta = (\zeta_1, \ldots, \zeta_n) \in \R^n$. A natural space for this operator is a weighted $L_1$-space:

\[
L_s:=\left\{u:\R^n\rightarrow\R:\int\limits_{\R^n}\frac{\vert u(x)\vert}{1+\vert x\vert^{n+2s}}dx<+\infty\right\}.
\]
The norm in $L_s$ is naturally given by
\[
\Vert u\Vert_{L_s}=\int\limits_{\R^n}\frac{\vert u(x)\vert}{1+\vert x\vert^{n+2s}}dx.
\]

The  study of system (\ref{1}) was mainly motivated  from the well known  fractional Lane-Emden problem

\begin{equation}\label{2}
\left\{
\begin{array}{rrll}
(-\Delta)^{s} u &=& u^p & {\rm in} \ \ \Omega\\
u &=& 0  & {\rm in} \ \ \R^n\setminus\Omega
\end{array}
\right.,
\end{equation}
where $\Omega$ is a smooth bounded open subset of $\R^{n}$, $n \geq 1$ and $0 < s < 1$.

Recently, it has been proved in \cite{SV} that this problem admits at least one positive solution for $1 < p < \frac{n + 2s}{n - 2s}$. The nonexistence has been established in \cite{ROS2} whenever $p \geq \frac{n + 2s}{n - 2s}$ and $\Omega$ is star-shaped. These results were known long before for $s = 1$, see the classical references \cite{AR, GS, ros237}.

For system of the type (\ref{1}) with $p=0=\theta$ and $r,q< 0$, existence results of positive solutions have been established when $qr>1$ in \cite{EM} for $s \neq t$ and in \cite{EM1} for $s = t$. The latter also proves existence and uniqueness of positive solution in the case that $ qr<1$. Finally, when $qr = 1$, the behavior of (\ref{1}) is resonant and the related eigenvalue problem has been studied in \cite{EM2}.

Nowadays, there has been some interest in systems of the type (\ref{1}) with $p,\theta\geq 0$ and $q, r > 0$. In \cite{ghergu}, the author studied existence, nonexistence, uniqueness, and regularity of solutions for the system (\ref{1}) with $s=1=t$.

In this paper, we are going to treat the system (\ref{1}) in the case $p,\theta\geq 0$ and $q, r > 0$. In this structure, the system above corresponds to the prototype equation (\ref{2}) in which the exponent $p$ is negative and generalize the results obtained in \cite{ghergu}. It is well known that for such a range of exponents, the system (\ref{1}) does not have a variational structure. To overcome this, we employ the sub-super method, which our approach relies on the boundary behavior of solutions to (\ref{2}) (with $p <0$) or more generally, to singular elliptic problems of the type
\begin{equation}\label{PS}
\left\{
\begin{array}{rrll}
(-\Delta)^{s} u &=& K(x)u^{-p},\ \ u>0  & {\rm in} \ \ \Omega\\
u &=& 0  & {\rm in} \ \ \mathbb{R}^n\setminus\Omega
\end{array}
\right.
\end{equation}
where $K\in C^{\nu}_{loc}(\Omega), \nu\in(0, 1)$, such that $\inf\limits_{\Omega} K>0$ and satisfies for some $0 \leq q<2s$ and $C_1, C_2>0$
\[
C_1d(x)^{-q} \leq K(x) \leq C_2d(x)^{-q}\text{ in }\Omega,
\]
where $d(x)=dist(x,\partial\Omega)$, studied by Adimurthi, Giacomoni and Santra in \cite{AGS}.

We say that a pair $(u,v)$ of continuous function in $\Omega$ and bounded in $\R^n$ is a positive classical
solution of system (\ref{1}), if $(-\Delta)^{s} u(x)$ and $(-\Delta)^{t} v(x)$ are well defined for all $x\in\Omega$, further $u$ and $v$ are positive in $\Omega$ and all equalities in (\ref{1}) hold pointwise in each corresponding set. Positive classical super and subsolutions are defined similarly.

We will establish our first result concerning the system (\ref{1}).

\begin{teor}\label{TC2} (Nonexistence). Let $p,\theta\geq 0$, $r,q>0.$ Then  the system (\ref{1}) has no positive classical solutions,  provided
	 that one of the following conditions holds:
\begin{itemize}
\item [$(i)$] $\frac{qt}{s}+p<1$ and $r\geq\frac{2t}{s}$;
\item [$(ii)$] $\frac{qt}{s}+p>1$ and $r(2s-qt)\geq 2s(1+p)$;
\item [$(iii)$] $\frac{rs}{t}+\theta<1$ and $q\geq\frac{2s}{t}$;
\item [$(iv)$] $\frac{rs}{t}+\theta>1$ and $q(2t-rs)\geq 2t(1+\theta)$;
\item [$(v)$] $p>\max\{1,\frac{rs}{t}-1\},\frac{2rs}{t}>(1-\theta)(1+p)$ and $qt(1+p-\frac{rs}{t})>(1+p)(1+\theta)s$;
\item [$(vi)$] $\theta>\max\{1,\frac{qt}{s}-1\},\frac{2qt}{s}>(1-p)(1+\theta)$ and $rs(1+\theta-\frac{qt}{s})>(1+p)(1+\theta)t$.
\end{itemize}

\end{teor}

\begin{rem}
The conditions (i) and (ii) in Theorem \ref{TC2} impose conditions on  the exponent $q$ to vary on the interval $(0, \frac{2s}{t})$, while in (v) the exponent $q$ can take any value greater than $\frac{2s}{t}$, provided  adjusting the other three exponents $p, r,\theta$ conveniently.  Finally, from  the conditions (iii), (iv) and (vi), the exponent $r$ is also restricted as above.

\end{rem}

Define the following quantities

\[
\alpha = p + \frac{qt}{s}\min\left\{1,\frac{2t -rs}{(1+ \theta)t}\right\},\ \ \ \beta= \theta +\frac{rs}{t}\min\left\{1,\frac{2s- qt}{(1+p)s}\right\}.
\]
These above  quantities  $\alpha$ and $\beta$ are related to the boundary behavior of the solution to the singular elliptic problem (\ref{PS}),  as they will be  explained in Proposition \ref{P2.6} below.

Next, we will state the existence of classical solutions to (\ref{1}).

\begin{teor}\label{TC1} (Existence). Let $p,\theta\geq 0$, $q,r>0$ satisfying the inequality
	
	\begin{equation}\label{5}
	(1 + p)(1+ \theta)- qr >0.
	\end{equation}
	
	 In addition, assume that  one of the following conditions below holds:

\begin{itemize}
\item [$(i)$] $\alpha\leq 1$ and $r<\frac{2t}{s}$;
\item [$(ii)$] $\beta\leq 1$ and $q<\frac{2s}{t}$;
\item [$(iii)$] $p,\theta\geq 1$, $r<\frac{2t}{s}$ and $q<\frac{2s}{t}$.
\end{itemize}
Then, the system (\ref{1}) has at least one positive classical solution $(u,v)\in (C^\eta(\R^n))^2$, for some $\eta\in(0,1)$.
\end{teor}

The proof  is made invoking  the Schauder's fixed point theorem in a suitable chosen closed convex subset of $(C^\eta(\R^n))^2$, for some $\eta\in(0,1)$, which contains all the functions having a certain rate of decay expressed in terms of the distance function $d(x)$ up to the boundary of $\Omega$.

The following necessary and sufficient
conditions for the existence of classical solutions to (\ref{1}) follows directly from Theorem \ref{TC2}(i) and (iii) and Theorem \ref{TC1}(i) and (ii).

\begin{cor} Let $p,\theta\geq 0$, $q, r >0$ satisfy (\ref{5}).
\begin{itemize}
\item [$(i)$] Assume $\frac{qt}{s}+p<1$. Then system (\ref{1}) has positive classical solutions if and only if $r<\frac{2t}{s}$;
\item [$(ii)$] Assume $\frac{rs}{t}+\theta<1$. Then system (\ref{1}) has positive classical solutions if and only if $q<\frac{2s}{t}$.
\end{itemize}
\end{cor}

\begin{teor}\label{TC3} (Uniqueness). Let $p,\theta\geq 0$, $q,r>0$, satisfy (\ref{5}) and one of the following conditions:
\begin{itemize}
\item [$(i)$] $\frac{qt}{s}+p<1$ and $r<\frac{2t}{s}$;
\item [$(ii)$] $\frac{rs}{t}+\theta<1$ and $q<\frac{2s}{t}$.
\end{itemize}
Then, the system (\ref{1}) has a unique positive classical solution.
\end{teor}

Several methods have been employed in the proof of existence, nonexistence and uniqueness results of positive solutions of elliptic systems. Our approach is inspired by a method developed by Ghergu in \cite{ghergu} to treat systems involving Laplace operators based on boundary behavior of the solution to (\ref{PS}), when $s=1$. Particularly, the boundary behavior of the solution to (\ref{PS}), proved by Adimurthi, Giacomoni, and Santra \cite{AGS}, as well as some fundamental results to be proved in the next section will play an important role in the proofs of Theorems of this work.

The paper is organized as follows. In Section 2 we obtain some preliminary properties related to the boundary behavior of the solution to (\ref{PS}). The rest of the Sections  are devoted to the proofs of our results.

\section{Notation and auxiliary results }
Consider the nonlocal eigenvalue problem

\begin{equation}\label{1.1}
\left\{
\begin{array}{rrll}
(-\Delta)^{s} u &=& \lambda u & {\rm in} \ \ \Omega\\
u &=& 0  & {\rm in} \ \ \R^n\setminus\Omega
\end{array}
\right..
\end{equation}
Since the operator $(-\Delta)^{s}$ is self-adjoint, by using a weak formulation and a suitable variational framework, Servadei and Valdinoci \cite{ros265} investigated in detail the discrete spectrum of $(-\Delta)^{s}$ in $\Omega$ for any $s \in (0,1)$. In particular, they proved that the first eigenvalue $\lambda_{1}(s):= \lambda_1((-\Delta)^{s})$ is positive, simple and characterized by
\[
\lambda_1(s)=\inf_{u\in X(\Omega)\setminus\{0\}}\frac{\int\limits_{\R^n}\vert(-\Delta)^{\frac{s}{2}} u\vert^2dx }{\int\limits_{\R^n}\vert u\vert^2}dx,
\]
where
\[
X(\Omega):=\{u\in H^s(\R^n):u=0 \text{ a.e in }\R^n\setminus\Omega\}.
\]

Let $\varphi_s$ be a nonnegative eigenfunction corresponding to $\lambda_{1}(s)$ in the weak sense. Results of H\"{o}lder regularity to the operator $(-\Delta)^{s}$ obtained by Ros-Oton and Serra \cite{ROS} imply that $\varphi_s \in C^{s}(\R^n)$ and moreover is a classical solution of \eqref{1.1} which is positive in $\Omega$. The last claim follows from Silvestre's strong maximum principle \cite{S} which holds for classical supersolutions (subsolutions).

By suitable normalization we may assume $|\varphi_s|_{\infty}=1$. In addition, it follows from the results in \cite{ROS1} that

\begin{equation}\label{8}cd(x)^s\leq \varphi_s(x)\leq\frac{1}{c}d(x)^s,\end{equation}
for some positive constant $c$.

We denote by $G_s(\cdot,\cdot)$ the Green's function of the fractional Laplace operator $(-\Delta)^{s}$ on $\Omega$. Let $w$ be a weak solution of the following problem
\begin{equation}\label{w}
\left\{
\begin{array}{rrll}
(-\Delta)^{s} w &=& h & {\rm in} \ \ \Omega\\
w &=& 0  & {\rm in} \ \ \R^n\setminus\Omega
\end{array}
\right..
\end{equation}

If $h \in C^\alpha_{loc}(\Omega)$, for some $\alpha\in (0,1)$, by Theorem 2.5 of \cite{xia}, there exists $\gamma>0$ such that $w\in C^{2s+\gamma}_{loc}(\Omega)$ is a classical solution of (\ref{w}), i.e, both equalities hold pointwise in each corresponding set. Therefore,
\begin{equation}\label{w2}
w(x)=\int\limits_{\Omega}G_s(x,y)h(y)dy\text{ in }\Omega\text{ and }w(x)=0\text{ in }\R^n\setminus\Omega.	
\end{equation}

Reciprocally, if $h \in C^\alpha(\overline{\Omega})$, for some $\alpha\in (0,1)$, by Theorem 1.2.3 of \cite{ros1} the function defined by setting (\ref{w2}) belongs to $C^{2s+\varepsilon}_{loc}(\Omega)\cap C(\overline{\Omega})\cap L_s$, fulfills $d(x)^{1-s}w\in C(\overline{\Omega})$, and $w$ is the only classical solution of problem (\ref{w}).

Now, let $\phi_s$  be the function that satisfies
\[
\left\{
\begin{array}{rrll}
(-\Delta)^{s} \phi_s &=& 1 & {\rm in} \ \ \Omega\\
\phi_s &=& 0  & {\rm in} \ \ \R^n\setminus\Omega
\end{array}
\right..
\]

By Silvestre's strong maximum principle (see \cite{S}), we get $\phi_s(x)>0$ in $\Omega$. Therefore,
\[\varphi_s(x)=\lambda_1(s)\int\limits_{\Omega}G_s(x,y)\varphi_s(y)dy\]
and
\[
\phi_s(x)=\int\limits_{\Omega}G_s(x,y)dy,
\]
which, as a consequence of the normalization of $\varphi_s$, leads to
 \begin{equation}\label{w3}
\varphi_s\leq \lambda_1(s)\phi_s.	
\end{equation}

An important tool for the uniqueness result of solutions of the system (\ref{1}) is as follows:

\begin{propo}\label{P2.1} Let $p\geq 0$ and $\psi:\Omega\rightarrow(0,\infty)$ be a continuous function. If $\underline{u}$ is a positive classical subsolution and $\overline{u}$ is a positive classical supersolution of
\[
\left\{
\begin{array}{rrll}
(-\Delta)^{s} u &=& \psi(x)u^{-p} & {\rm in} \ \ \Omega\\
u &=& 0  & {\rm in} \ \ \mathbb{R}^n\setminus\Omega
\end{array}
\right.,
\]
then $\underline{u}\leq\overline{u}$ in $\Omega$.
\end{propo}

\n {\bf Proof.} If $p=0$ the result is a consequence of the  Silvestre's strong maximum principle. Suppose  $p >0,$ and assume by contradiction that the set $\omega := \{x\in\Omega:\overline{u}(x) < \underline{u}(x)\}$ is not empty and let $w := \underline{u} -\overline{u}$. Then, $w$ achieves its maximum on $\Omega$ at a point $x_0\in \omega$. Then,
\begin{eqnarray*}
0&\leq& C(n,s)\lim_{\varepsilon\searrow 0}\int\limits_{\R^{n}\setminus B_\varepsilon(x_0)}\frac{w(x_0)-w(y)}{\vert x_0-y\vert^{n+2s}}\; dy=(-\Delta)^{s}w(x_0)\\ &\leq&\psi(x_0)[\underline{u}(x_0)^{-p}-\overline{u}(x_0)^{-p}]<0,
\end{eqnarray*}
which is a contradiction. Therefore, $\omega=\emptyset$, that is, $\underline{u}\leq\overline{u}$ in $\Omega$. \fim\\

Now an important tool for the nonexistence and uniqueness results of solutions of the system (\ref{1}) is as follows:

\begin{propo}\label{P2.3}
Let $(u, v)$ be a positive classical solution of system (\ref{1}). Then, there exists a constant $c > 0$ such
that
\begin{equation}\label{6}
u(x)\geq cd(x)^s\text{ and }v(x) \geq cd(x)^t \text{ in }\Omega.
\end{equation}
\end{propo}
\n {\bf Proof.}
Let $(u, v)$ be a positive classical solution of \eqref{1}. By inequalities (\ref{8}) and (\ref{w3}), there is a constant $c_0>0$ such that $\phi_s(x)\geq c_0 d(x)^s$ and $\phi_t(x)\geq c_0 d(x)^t$ in $\Omega$. Notice that $(-\Delta)^{s}u\geq C =(-\Delta)^{s}(C\phi_s)$ in $\Omega$, where $C=\min\limits_{\Omega}\{u^{-p}v^{-q}\}>0$. Then, by Silvestre's strong maximum principle, we deduce $u(x)\geq C\phi_s(x) \geq cd(x)^s$ in $\Omega$ and similarly $v(x)\geq cd(x)^t$ in $\Omega$, where $c>0$ is a positive constant. \fim \\

The following result is a direct consequence of Silvestre's strong maximum principle, inequality (\ref{8}) and Theorem 1.2 of \cite{AGS}. This is the key tool for the existence, nonexistence and uniqueness results of solutions of the system (\ref{1}).

\begin{propo}\label{P2.6}
Let $p\geq 0$ and $\gamma>0$. There are constants $c, C>0$ such that any positive classical subsolution $\underline{u}$ and any positive classical supersolution $\overline{u}$ of problem
\begin{equation}\label{13}
\left\{
\begin{array}{rrll}
(-\Delta)^{s} u &=& d(x)^{-\gamma}u^{-p} & {\rm in} \ \ \Omega\\
u &=& 0  & {\rm in} \ \ \mathbb{R}^n\setminus\Omega
\end{array}
\right.,
\end{equation}
satisfies:
\begin{itemize}
\item [$(i)$] $\underline{u}(x)\leq Cd(x)^s$ and $\overline{u}\geq cd(x)^s$ in $\Omega$, if $\frac{\gamma}{s}+p<1$;
\item [$(ii)$] $\underline{u}(x)\leq Cd(x)^s\left(\ln\left(\frac{2}{\varphi_s}\right)\right)^{\frac{1}{1+p}}$ and $\overline{u}\geq cd(x)^s\left(\ln\left(\frac{2}{\varphi_s}\right)\right)^{\frac{1}{1+p}}$ in $\Omega$, if $\frac{\gamma}{s}+p=1$;
\item [$(iii)$] $\underline{u}(x)\leq Cd(x)^{\frac{2s-\gamma}{1+p}}$ and $\overline{u}\geq cd(x)^{\frac{2s-\gamma}{1+p}}$ in $\Omega$, if $\frac{\gamma}{s}+p>1$ with $0<\gamma<2s$.
\end{itemize}
\end{propo}

Finally, Theorem 1.2(iii) of \cite{AGS} also guarantees that the problem (\ref{13}) has no positive classical solution, if $\gamma\geq 2s$. Such claimed is important for the proof of nonexistence results of positive classical solutions of the system (\ref{1}).

\section{Proof of Theorem \ref{TC2}}

Notice that the system (\ref{1}) is invariant under the transform $(u, v,p, q, r, \theta, s)\rightarrow(v,u, \theta, r, q,p, t)$, so that,  we  need to prove only the cases {\bf (i), (ii)} and {\bf (v)}.

Suppose  that there exists $(u,v)$ a positive classical solution of system (\ref{1}). By Proposition \ref{P2.3}, we can find $c>0$ such that (\ref{6}) holds.\\

{\bf (i)} $\frac{qt}{s}+p<1$ and $r\geq \frac{2t}{s}$. Using the estimate (\ref{6}) in the first equation of the system (\ref{1}) we have

\begin{equation}\label{22}
\left\{
\begin{array}{rrll}
(-\Delta)^{s} u &=& v^{-q}u^{-p}\leq c_1d(x)^{-qt}u^{-p},\ \ \ u>0 & {\rm in} \ \ \Omega\\
u &=& 0  & {\rm in} \ \ \mathbb{R}^n\setminus\Omega
\end{array}
\right.,
\end{equation}
for some $c_1>0$. By Proposition \ref{P2.6}(i) we conclude  $u(x)\leq c_2d(x)^s$ in $\Omega$, for some $c_2>0$. From this and (\ref{6}), we have there exists $c_0,c_3>0$ such that $c_0d(x)^{-rs}\leq u^{-r}\leq c_3d(x)^{-rs}$ in $\Omega$. Using the second equation of (\ref{1}) we find

\begin{equation}\label{23}
\left\{
\begin{array}{rrll}
(-\Delta)^{t} v &=&u^{-r}v^{-\theta},\ \ \ v>0 & {\rm in} \ \ \Omega\\
v &=& 0  & {\rm in} \ \ \mathbb{R}^n\setminus\Omega
\end{array}
\right..
\end{equation}
According to Theorem 1.2(iii) of \cite{AGS}, this is impossible, since $rs\geq 2t$.\\

{\bf (ii)} $\frac{qt}{s}+p>1$ and $r(2s-qt)\geq 2t(1+p)$. In the same manner as above, $u$ satisfies the problem (\ref{22}). From Proposition \ref{P2.6}(iii), we now deduce
\[
u(x)\leq c_1d(x)^{\frac{2s-qt}{1+p}}
\]
in $\Omega$, for some $c_1>0$. Since $\frac{qt}{s}+p>1$, we have $\frac{2s-qt}{1+p}>s$. From this and (\ref{6}) we deduce
\[
u(x)\geq cd(x)^{\frac{2s-qt}{1+p}}.
\]

Then, there are $c_2,c_3>0$ such that
\[
c_2d(x)^{-\frac{r(2s-qt)}{1+p}}\leq u^{-r}\leq c_3d(x)^{-\frac{r(2s-qt)}{1+p}}\text{ in }\Omega.
\]
Now, using the second equation of (\ref{1}) we have $v$ is a classical solution of problem (\ref{23}), which is impossible in view of Theorem 1.2(iii) of \cite{AGS}, since $\frac{r(2s-qt)}{1+p}\geq 2t$.\\

{\bf (v)} Let $M=\sup\limits_{x\in\overline{\Omega}}v$. From the first equation of the system (\ref{1}) we find

\[
\left\{
\begin{array}{rrll}
(-\Delta)^{s} u &\geq& M^{-q}u^{-p},\ \ \ u>0 & {\rm in} \ \ \Omega\\
u &=& 0  & {\rm in} \ \ \mathbb{R}^n\setminus\Omega
\end{array}
\right..
\]

From Proposition \ref{P2.6}(iii), we have $u(x)\geq c_1d(x)^{\frac{2s}{1+p}}$ in $\Omega$, for some $c_1>0$. Combining this estimate with the second equation of (\ref{1}) we have

\[
\left\{
\begin{array}{rrll}
(-\Delta)^{t} v &\leq& c_2 d(x)^{-\frac{2rs}{1+p}}v^{-\theta},\ \ \ v>0 & {\rm in} \ \ \Omega\\
v &=& 0  & {\rm in} \ \ \mathbb{R}^n\setminus\Omega
\end{array}
\right.,
\]
where $c_2>0$. Since $\frac{2rs}{t(1+p)}+\theta>1$, again by Proposition \ref{P2.6}(iii) we obtain that the function $v$ satisfies
\[
v(x)\leq c_3d(x)^{\frac{2t(1+p)-2rs}{(1+p)(1+\theta)}}\text{ in }\Omega,
\]
for some $c_3>0$. Since $\frac{2rs}{t}>(1-\theta)(1+p)$, we have $\frac{2t(1+p)-2rs}{(1+p)(1+\theta)}>t$. From this and (\ref{6}) we deduce
\[
u(x)\geq cd(x)^{\frac{2t(1+p)-2rs}{(1+p)(1+\theta)}}.
\]

Then, there exists $c_4,c_5>0$ such that
\[
c_4d(x)^{-\frac{q(2t(1+p)-2rs)}{(1+p)(1+\theta)}}\leq v^{-q}\leq c_5d(x)^{-\frac{q(2t(1+p)-2rs)}{(1+p)(1+\theta)}}.
\]
Now, using the first equation of (\ref{1}) we have $u$ is a classical solution of problem
\[
\left\{
\begin{array}{rrll}
(-\Delta)^{s} u &=& v^{-q}u^{-p},\ \ \ u>0 & {\rm in} \ \ \Omega\\
u &=& 0  & {\rm in} \ \ \mathbb{R}^n\setminus\Omega
\end{array}
\right.,
\]
which contradicts Theorem 1.2(iii) of \cite{AGS}, since $qt(1+p-\frac{rs}{t})>(1+p)(1+\theta)s$. Thus, the system (\ref{1}) has no positive classical solutions. This completes the proof of Theorem \ref{TC2}. \fim \\

\section{Proof of Theorem \ref{TC1}}

{\bf (i)} The proof is made in six cases according to bounded behavior of singular elliptic problems of the type (\ref{PS}), as it was pointed out in Proposition \ref{P2.6}.\\

{\bf Case 1:} $\frac{rs}{t} +\theta > 1$ and $\alpha =\frac{q(2t-rs)}{s(1+\theta)}+p < 1$. From Proposition \ref{P2.6}(i) and (iii) there exist $0 < c_1 <1 < c_2$ such that:\\

$\bullet$ Any positive classical subsolution $\underline{u}$ and any positive classical supersolution $\overline{u}$ of the problem
\begin{equation}\label{27}
\left\{
\begin{array}{rrll}
(-\Delta)^{s} u &=& d(x)^{-\frac{q(2t-rs)}{(1+\theta)}}u^{-p} & {\rm in} \ \ \Omega\\
u &=& 0  & {\rm in} \ \ \mathbb{R}^n\setminus\Omega
\end{array}
\right.,
\end{equation}
satisfy
\begin{equation}\label{28}
\overline{u}(x) \geq c_1d(x)^s\text{ and }\underline{u}(x) \leq c_2d(x)^s\text{ in }\Omega. \\
\end{equation}

$\bullet$ Any positive classical subsolution $\underline{v}$ and any positive classical supersolution $\overline{v}$ of the problem
\begin{equation}\label{29}
\left\{
\begin{array}{rrll}
(-\Delta)^{t} v &=& d(x)^{-rs}v^{-\theta} & {\rm in} \ \ \Omega\\
v &=& 0  & {\rm in} \ \ \mathbb{R}^n\setminus\Omega
\end{array}
\right.,
\end{equation}
satisfy
\[
\overline{v}(x)\geq c_1d(x)^{\frac{2t-rs}{1+\theta}}\text{ and }\underline{v}(x)\leq c_2 d(x)^{\frac{2t-rs}{1+\theta}}\text{ in }\Omega.
\]

We fix $0<m_1 < 1<M_1$ and $0<m_2 < 1<M_2$ such that
\begin{equation}\label{31}
M_1^{\frac{r}{1+\theta}} m_2 \leq c_1 < c_2 \leq M_1m_2^{\frac{q}{1+p}}
\end{equation}
and
\begin{equation}\label{32}
M_2^{\frac{q}{1+p}} m_1 \leq c_1 < c_2 \leq M_2m_1^{\frac{r}{1+\theta}}.
\end{equation}
Note that the above choice of $m_i,M_i$ $(i = 1, 2)$ is possible in view of (\ref{5}).

Let $\varepsilon_1 >0$ small enough. Here $X$ stands for the Banach space
\[
\{(u,v) \in C^{s-\varepsilon_1}(\R^n)\times C^{\frac{2t-rs}{1+\theta}-\varepsilon_1}(\R^n) : u = v = 0\ {\rm in}\ \R^n\setminus\Omega\}
\]
endowed with the product norm

\[
\| (u,v) \|_X := \| u \|_{C^{s-\varepsilon_1}(\R^n)} + \| v \|_{C^{\frac{2t-rs}{1+\theta}-\varepsilon_1}(\R^n)}\, .
\]

Set
\begin{eqnarray*}
\mathcal{A}:=\left\{(u,v)\in X : \begin{array}{c}  m_1d(x)^s\leq u\leq M_1d(x)^s    \\  \text{ and } \\    m_2d(x)^{\frac{2t-rs}{1+\theta}}\leq v\leq M_2d(x)^{\frac{2t-rs}{1+\theta}}\text{ in }\Omega \end{array} \right\}.
\end{eqnarray*}

For any $(u, v) \in\mathcal{A}$, let  $(T u, T v)$ be  the unique positive classical solution of the decoupled system
\begin{equation}\label{33}
\left\{
\begin{array}{llll}
(-\Delta)^{s}(Tu) = v^{-q}(Tu)^{-p} & {\rm in} \ \ \Omega\\
(-\Delta)^{t}(Tv) = u^{-r}(Tv)^{-\theta} & {\rm in} \ \ \Omega\\
Tu= Tv=0 & {\rm in} \ \ \R^n\setminus\Omega
\end{array}
\right.
\end{equation}
and define
\begin{equation}\label{34}
\mathcal{F}:\mathcal{A}\rightarrow X \text{ by }\mathcal{F}(u, v) = (T u, T v)\text{ for any }(u, v)\in\mathcal{A}.
\end{equation}

It is proved in \cite{AGS}, the existence of positive classical solution $Tu\in C^{s}(\R^n)$ and $Tv\in C^{\frac{2t-rs}{1+\theta}}(\R^n),$ and the uniqueness of the positive weak solution in each equation of the system (\ref{33}). We define the space $X$ as subspace of
\[
C^{s-\varepsilon_1}(\R^n)\times C^{\frac{2t-rs}{1+\theta}-\varepsilon_1}(\R^n),
\]
for some $\varepsilon_1 >0$ small enough, to ensure the compactness of the operator $\mathcal{F}$ (see Step 2 below).

Therefore, if  $\mathcal{F}$ has a fixed point in $\mathcal{A}$, then the existence of a positive classical solution to system (\ref{1}) follows. To this end, we shall prove that $\mathcal{F}$ satisfies the conditions:
\[
\mathcal{F}(\mathcal{A}) \subseteq\mathcal{A},\ \ \ \mathcal{F}\text{ is compact and continuous}.
\]
Hence, by Schauder's fixed point theorem we deduce that $\mathcal{F}$ has a fixed point in $\mathcal{A}$, which is a positive classical solution to (\ref{1}).\\

{\bf Step 1:} $ \mathcal{F}(\mathcal{A}) \subseteq\mathcal{A}$. Take $(u, v)\in\mathcal{A}$. From the inequality
\[
v\leq M_2d(x)^{\frac{2t-rs}{1+\theta}}\text{ in }\Omega,
\]
we obtain that $Tu$ satisfies
\[
\left\{
\begin{array}{rrll}
(-\Delta)^{s} (Tu) &\geq& M_2^{-q}d(x)^{-\frac{q(2t-rs)}{1+\theta}}(Tu)^{-p},\ \ \ Tu>0 & {\rm in} \ \ \Omega\\
Tu &=& 0  & {\rm in} \ \ \mathbb{R}^n\setminus\Omega
\end{array}
\right..
\]
Thus, $\overline{u}:=M_2^{\frac{q}{1+p}}Tu$ is a positive classical supersolution to (\ref{27}), because $-q+\frac{q}{1+p}=-p\frac{q}{1+p}$. By (\ref{28}) and (\ref{32}) we obtain
\[
Tu =M_2^{-\frac{q}{1+p}}\overline{u}\geq c_1M_2^{-\frac{q}{1+p}} d(x)^s\geq m_1d(x)^s\text{ in }\Omega.
\]
By inequality $m_2d(x)^{\frac{2t-rs}{1+\theta}}\leq v$ in $\Omega$ and the definition of $Tu$ we conclude that
\[
\left\{
\begin{array}{rrll}
(-\Delta)^{s} (Tu) &\leq& m_2^{-q}d(x)^{-\frac{q(2t-rs)}{1+\theta}}(Tu)^{-p},\ \ \ Tu>0 & {\rm in} \ \ \Omega\\
Tu &=& 0  & {\rm in} \ \ \mathbb{R}^n\setminus\Omega
\end{array}
\right..
\]
Therefore, $\underline{u} := m_2^{\frac{q}{1+p}}Tu$ is a positive classical subsolution of problem (\ref{27}). Hence, from (\ref{28}) and (\ref{31}) we have
\[
Tu =m_2^{-\frac{q}{1+p}}\overline{u}\leq c_2m_2^{-\frac{q}{1+p}} d(x)^s\leq M_1d(x)^s\text{ in }\Omega.
\]
This way, we have proved that $Tu$ satisfies
\[
m_1d(x)^s\leq Tu\leq M_1d(x)^s\text{ in }\Omega.
\]
Similarly, using the definition of $\mathcal{A}$ and the properties of the sub and supersolutions of problem (\ref{29}) we can prove that $Tv$ satisfies
\[
m_2d(x)^{\frac{2t-rs}{1+\theta}}\leq Tv\leq M_2d(x)^{\frac{2t-rs}{1+\theta}}\text{ in }\Omega.
\]
Then, $(Tu, Tv)\in\mathcal{A}$ for all $(u, v)\in\mathcal{A}$, that is, $\mathcal{F}(\mathcal{A}) \subseteq\mathcal{A}$.\\

{\bf Step 2:} $\mathcal{F}$ is compact and continuous. Let $(u, v)\in\mathcal{A}$. Then, we conclude $Tu\in C^{s}(\R^n)$ and $Tv\in C^{\frac{2t-rs}{1+\theta}}(\R^n)$. Recalling that the embedding $C^{0,s}(\overline{\Omega})\hookrightarrow C^{0,s-\varepsilon_1}(\overline{\Omega})$ and $C^{0,\frac{2t-rs}{1+\theta}}(\overline{\Omega})\hookrightarrow C^{0,\frac{2t-rs}{1+\theta}-\varepsilon_1}(\overline{\Omega})$ are compact, it follows that $\mathcal{F}$ is also compact.\\

Now,rest to prove that $\mathcal{F}$ is continuous. To this end, let $(u_n, v_n)\subset\mathcal{A}$ be such that $u_n\rightarrow u$ in $C^{s-\varepsilon_1}(\R^n)$ and $v_n\rightarrow v$ in $C^{\frac{2t-rs}{1+\theta}-\varepsilon_1}(\R^n)$ as $n\rightarrow\infty$. Since  $\mathcal{F}$ is compact, there exists $(U, V)\in\mathcal{A}$ such that up to a subsequence we get
\[
Tu_n\rightarrow U\text{ in }C^{s-\varepsilon_1}(\R^n)\text{ and }Tv_n\rightarrow V\text{ in }C^{\frac{2t-rs}{1+\theta}-\varepsilon_1}(\R^n)\text{ as }n\rightarrow\infty.
\]
By Theorem 2.7 of \cite{xia}, we have $(U, V)$ is a positive viscosity solution of system }(see definition in the paper \cite{xia}).
\[
\left\{
\begin{array}{llll}
(-\Delta)^{s}U = v^{-q}U^{-p} & {\rm in} \ \ \Omega\\
(-\Delta)^{t}V = u^{-r}V^{-\theta} & {\rm in} \ \ \Omega\\
U= V=0 & {\rm in} \ \ \R^n\setminus\Omega
\end{array}
\right..
\]
From the  uniqueness of positive weak solution of the problem (\ref{33}), it follows that $Tu=U$ and $Tv=V$. So,
\[
Tu_n\rightarrow Tu\text{ in }C^{s-\varepsilon_1}(\R^n)\text{ and }Tv_n\rightarrow Tv\text{ in }C^{\frac{2t-rs}{1+\theta}-\varepsilon_1}(\R^n)\text{ as }n\rightarrow\infty.
\]
So that,  $\mathcal{F}$ is continuous.

Applying the Schauder's fixed point theorem, there exists $(u, v)\in\mathcal{A}$ such that $\mathcal{F}(u, v)=(u, v)$, that is, $Tu = u$ and $Tv=v$. Therefore, $(u, v)$ is a positive classical solution of system (\ref{1}).\\

The others cases will be considered  similarly. But, due to the different boundary behavior of solutions described in Proposition \ref{P2.6}, the set $\mathcal{A}$ and the constants $c_1, c_2$ have to be modified accordingly. We shall point out how  these constants are chosen in order to apply the
Schauder's fixed point theorem.\\

{\bf Case 2:} $\frac{rs}{t} +\theta = 1$ and $\alpha =\frac{qt}{s}+p < 1$. By Proposition \ref{P2.6}(i) and (ii) there exists $0<a<1$ and $0 < c_1 <1 < c_2$ such that:\\

$\bullet$ Any positive classical subsolution $\underline{u}$ of the problem
\[
\left\{
\begin{array}{rrll}
(-\Delta)^{s} u &=& d(x)^{-qt}u^{-p} & {\rm in} \ \ \Omega\\
u &=& 0  & {\rm in} \ \ \mathbb{R}^n\setminus\Omega
\end{array}
\right.,
\]
verify
\[
\underline{u}(x) \leq c_2d(x)^s\text{ in }\Omega. \\
\]

$\bullet$ Any positive classical supersolution $\overline{u}$ of the problem
\[
\left\{
\begin{array}{rrll}
(-\Delta)^{s} u &=& d(x)^{-qt(t-at)}u^{-p} & {\rm in} \ \ \Omega\\
u &=& 0  & {\rm in} \ \ \mathbb{R}^n\setminus\Omega
\end{array}
\right.,
\]
satisfy
\[
\overline{u}(x) \geq c_1d(x)^s\text{ in }\Omega. \\
\]

$\bullet$ Any positive classical subsolution $\underline{v}$ and any positive classical supersolution $\overline{v}$ of problem (\ref{29}) satisfy
\[
\overline{v}(x)\geq c_1d(x)^{t}\text{ and }\underline{v}(x)\leq c_2d(x)^{t-at}\text{ in }\Omega.
\]

Let $\varepsilon_1 >0$ small enough. Here $X$ stands for the Banach space
\[
\{(u,v) \in C^{s-\varepsilon_1}(\R^n)\times C^{t-\varepsilon-\varepsilon_1}(\R^n) : u = v = 0\ {\rm in}\ \R^n\setminus\Omega\},
\]
for any $\varepsilon >0$ small enough, endowed with the product norm

\[
\| (u,v) \|_X := \| u \|_{C^{s-\varepsilon_1}(\R^n)} + \| v \|_{C^{t-\varepsilon-\varepsilon_1}(\R^n)}\, .
\]

Set
\begin{eqnarray*}
\mathcal{A}:=\left\{(u,v)\in X :m_1d(x)^s\leq u\leq M_1d(x)^s\text{ and } m_2d(x)^{t}\leq v\leq M_2d(x)^{t-at}\text{ in }\Omega \right\},
\end{eqnarray*}
where $0<m_i < 1<M_i\ (i = 1, 2)$ satisfy (\ref{31}), (\ref{32}) and
\[
m_2[diam(\Omega)]^{at}<M_2.
\]
Define the operator $\mathcal{F}$  as in the Case 1 by (\ref{33}) and (\ref{34}). The inclusion $\mathcal{F}(\mathcal{A}) \subseteq\mathcal{A}$ and that $\mathcal{F}$ is continuous and compact follow as before.\\

{\bf Case 3:} $\frac{rs}{t} +\theta <1$ and $\alpha=\frac{qt}{s}+p <1$. Let $\varepsilon_1 >0$ small enough. Here $X$ stands for the Banach space
\[
\{(u,v) \in C^{s-\varepsilon_1}(\R^n)\times C^{t-\varepsilon_1}(\R^n) : u = v = 0\ {\rm in}\ \R^n\setminus\Omega\},
\]
endowed with the product norm

\[
\| (u,v) \|_X := \| u \|_{C^{s-\varepsilon_1}(\R^n)} + \| v \|_{C^{t-\varepsilon_1}(\R^n)}\, .
\]

In the same manner we define
\begin{eqnarray*}
\mathcal{A}:=\left\{(u,v)\in X :m_1d(x)^s\leq u\leq M_1d(x)^s\text{ and } m_2d(x)^{t}\leq v\leq M_2d(x)^{t}\text{ in }\Omega \right\},
\end{eqnarray*}
where $0<m_i < 1<M_i\ (i = 1, 2)$ satisfy (\ref{31}), (\ref{32}) for suitable constants $c_1$ and $c_2$.\\

{\bf Case 4:} $\frac{rs}{t} +\theta <1$ and $\alpha=\frac{qt}{s}+p =1$. The approach is the same as in Case 2 above.

Let $\varepsilon_1 >0$ small enough. Here $X$ stands for the Banach space
\[
\{(u,v) \in C^{s-\varepsilon-\varepsilon_1}(\R^n)\times C^{t-\varepsilon_1}(\R^n) : u = v = 0\ {\rm in}\ \R^n\setminus\Omega\},
\]
for any $\varepsilon >0$ small enough, endowed with the product norm

\[
\| (u,v) \|_X := \| u \|_{C^{s-\varepsilon-\varepsilon_1}(\R^n)} + \| v \|_{C^{t-\varepsilon_1}(\R^n)}\, .
\]

Set
\begin{eqnarray*}
\mathcal{A}:=\left\{(u,v)\in X :m_1d(x)^s\leq u\leq M_1d(x)^{s-as}\text{ and } m_2d(x)^{t}\leq v\leq M_2d(x)^{t}\text{ in }\Omega \right\},
\end{eqnarray*}
for some $0<a<1$, where $0<m_i < 1<M_i\ (i = 1, 2)$ satisfy (\ref{31}), (\ref{32}) and
\[
m_1[diam(\Omega)]^{as}<M_1.
\]

{\bf Case 5:} $\frac{rs}{t} +\theta >1$ and $\alpha=\frac{qt}{s}+p =1$. Let $0 < a < 1$ be fixed such that $\frac{ars}{t} +\theta >1$. Then,
\[
\frac{q(2t-rs)}{s(1+\theta)}+p<1\text{ and }\frac{q(2t-rsa)}{s(1+\theta)}+p<1.
\]

So, by Proposition \ref{P2.6}(i), (iii), there exist $0 < c_1 < 1 < c_2$ such that:\\

$\bullet$ Any positive classical subsolution $\underline{u}$ of the problem
\[
\left\{
\begin{array}{rrll}
(-\Delta)^{s} u &=& d(x)^{-\frac{q(2t-rsa)}{1+\theta}}u^{-p} & {\rm in} \ \ \Omega\\
u &=& 0  & {\rm in} \ \ \mathbb{R}^n\setminus\Omega
\end{array}
\right.,
\]
verify
\[
\underline{u}(x) \leq c_2d(x)^{sa}\text{ in }\Omega. \\
\]

$\bullet$ Any positive classical supersolution $\overline{u}$ of the problem
\[
\left\{
\begin{array}{rrll}
(-\Delta)^{s} u &=& d(x)^{-\frac{q(2t-rs)}{1+\theta}}u^{-p} & {\rm in} \ \ \Omega\\
u &=& 0  & {\rm in} \ \ \mathbb{R}^n\setminus\Omega
\end{array}
\right.,
\]
satisfy
\[
\overline{u}(x) \geq c_1d(x)^s\text{ in }\Omega. \\
\]

$\bullet$ Any positive classical subsolution $\underline{v}$ of problem (\ref{29}) satisfies
\[
\underline{v}(x)\leq c_2d(x)^{\frac{2t-rs}{1+\theta}}\text{ in }\Omega.
\]

$\bullet$ Any positive classical supersolution $\overline{v}$ of problem
\[
\left\{
\begin{array}{rrll}
(-\Delta)^{t} v &=& d(x)^{-ars}v^{-\theta} & {\rm in} \ \ \Omega\\
v &=& 0  & {\rm in} \ \ \mathbb{R}^n\setminus\Omega
\end{array}
\right.,
\]
satisfies
\[
\overline{v}(x)\geq c_1d(x)^{\frac{2t-rsa}{1+\theta}}\text{ in }\Omega.
\]
Let $\varepsilon_1 >0$ small enough. Here $X$ stands for the Banach space
\[
\{(u,v) \in C^{s-\varepsilon-\varepsilon_1}(\R^n)\times C^{\frac{2t-rs}{1+\theta}-\varepsilon_1}(\R^n) : u = v = 0\ {\rm in}\ \R^n\setminus\Omega\},
\]
for any $\varepsilon >0$ small enough, endowed with the product norm

\[
\| (u,v) \|_X := \| u \|_{C^{s-\varepsilon-\varepsilon_1}(\R^n)} + \| v \|_{C^{\frac{2t-rs}{1+\theta}-\varepsilon_1}(\R^n)}\, .
\]

Set
\begin{eqnarray*}
\mathcal{A}:=\left\{(u,v)\in X :m_1d(x)^s\leq u\leq M_1d(x)^{sa}\text{ and } m_2d(x)^{\frac{2t-rsa}{1+\theta}}\leq v\leq M_2d(x)^{\frac{2t-rs}{1+\theta}}\text{ in }\Omega \right\},
\end{eqnarray*}
where $0<m_i < 1<M_i\ (i = 1, 2)$ satisfy (\ref{31}), (\ref{32}) in which the constants $c_1, c_2$ are those given above and
\[
m_1[diam(\Omega)]^{s-as}<M_1\text{ and }m_2[diam(\Omega)]^{\frac{r(s-as)}{1+\theta}}<M_2.
\]

{\bf Case 6:} $\frac{rs}{t} +\theta =1$ and $\alpha=\frac{qt}{s}+p =1$. We proceed in the same manner as above by considering $X$ stands for the Banach space
\[
\{(u,v) \in C^{s-\varepsilon-\varepsilon_1}(\R^n)\times C^{t-\varepsilon-\varepsilon_1}(\R^n) : u = v = 0\ {\rm in}\ \R^n\setminus\Omega\},
\]
for any $\varepsilon, \varepsilon_1 >0$ small enough, endowed with the product norm

\[
\| (u,v) \|_X := \| u \|_{C^{s-\varepsilon-\varepsilon_1}(\R^n)} + \| v \|_{C^{t-\varepsilon-\varepsilon_1}(\R^n)}\, .
\]

Set
\begin{eqnarray*}
\mathcal{A}:=\left\{(u,v)\in X :m_1d(x)^s\leq u\leq M_1d(x)^{s-a_1s}\text{ and } m_2d(x)^{t}\leq v\leq M_2d(x)^{t-a_2t}\text{ in }\Omega \right\},
\end{eqnarray*}
where $0<a_1,a_2<1$ are fixed constants and $0<m_i < 1<M_i\ (i = 1, 2)$ satisfy (\ref{31}), (\ref{32}) for suitable constants $c_1, c_2>0$ and
\[
m_1[diam(\Omega)]^{a_1s}<M_1\text{ and }m_2[diam(\Omega)]^{a_2t}<M_2.
\]

{\bf (iii)} Let
\[
a=\frac{2st\left(\frac{1+\theta}{t}-\frac{q}{s}\right)}{(1+p)(1+\theta)-qr}\text{ and }b=\frac{2st\left(\frac{1+p}{s}-\frac{r}{t}\right)}{(1+p)(1+\theta)-qr}.
\]

Then
\begin{equation}\label{36}
a=\frac{2s-bq}{1+p}\text{ and }b=\frac{2t-ra}{1+\theta}.
\end{equation}
From hypothesis, we have $a<s$ and $b<t$. Then, $bq<qt<2s$ and $ar<rs<2t$. Now, since $\frac{bq}{s}+ p > 1$ and $\frac{ar}{t}+\theta > 1$, from Proposition \ref{P2.6}(iii) and (\ref{36}) above we can find $0 < c_1 <1 < c_2$ such that:\\

$\bullet$ Any positive classical subsolution $\underline{u}$ and any positive classical supersolution $\overline{u}$ of the problem
\[
\left\{
\begin{array}{rrll}
(-\Delta)^{s} u &=& d(x)^{-bq}u^{-p} & {\rm in} \ \ \Omega\\
u &=& 0  & {\rm in} \ \ \mathbb{R}^n\setminus\Omega
\end{array}
\right.,
\]
satisfy
\[
\overline{u}(x) \geq c_1d(x)^a\text{ and }\underline{u}(x) \leq c_2d(x)^a\text{ in }\Omega. \\
\]

$\bullet$ Any positive classical subsolution $\underline{v}$ and any positive classical supersolution $\overline{v}$ of the problem
\[
\left\{
\begin{array}{rrll}
(-\Delta)^{t} v &=& d(x)^{-ar}v^{-\theta} & {\rm in} \ \ \Omega\\
v &=& 0  & {\rm in} \ \ \mathbb{R}^n\setminus\Omega
\end{array}
\right.,
\]
verify
\[
\overline{v}(x)\geq c_1d(x)^{b}\text{ and }\underline{v}(x)\leq c_2d(x)^{b}\text{ in }\Omega.
\]

As before, let $\varepsilon_1 >0$ small enough and define $X$ to be the Banach space
\[
\{(u,v) \in C^{a-\varepsilon_1}(\R^n)\times C^{b-\varepsilon_1}(\R^n) : u = v = 0\ {\rm in}\ \R^n\setminus\Omega\}
\]
endowed with the product norm

\[
\| (u,v) \|_X := \| u \|_{C^{a-\varepsilon_1}(\R^n)} + \| v \|_{C^{b-\varepsilon_1}(\R^n)}\, .
\]

Set
\begin{eqnarray*}
\mathcal{A}:=\left\{(u,v)\in X :m_1d(x)^a\leq u\leq M_1d(x)^a\text{ and } m_2d(x)^{b}\leq v\leq M_2d(x)^{b}\text{ in }\Omega \right\},
\end{eqnarray*}
where $0<m_1 < 1 <M_1$ and $0<m_2 < 1 <M_2$ satisfy (\ref{31}) and (\ref{32}). This completes the proof of Theorem \ref{TC1}. \fim \\

\section{Proof of Theorem \ref{TC3}}

We shall prove only {\bf (i)}; the case {\bf (ii)} follows similarly.

Let $(u_1,v_1)$ and $(u_2,v_2)$ be two positive classical solutions of system (\ref{1}). Note that if $\frac{qt}{s}+p<1$, then by Theorem \ref{TC2}, we deduce $rs<2t$. By Proposition \ref{P2.3} there exists $c_1>0$ such that

\begin{equation}\label{410}
u_i\geq c_1d(x)^s\text{ and }v_i\geq c_1d(x)^t
\end{equation}
in $\Omega$, $i=1,2$. Then, $u_i$ satisfies

\[
\left\{
\begin{array}{rrll}
(-\Delta)^{s} u_i &=& v_i^{-q}u_i^{-p} \leq c_2 d(x)^{-qt}u_i^{-p},\ u_i>0 & {\rm in} \ \ \Omega\\
u_i &=& 0  & {\rm in} \ \ \mathbb{R}^n\setminus\Omega
\end{array}
\right.,
\]
for some $c_2>0$. Since $\frac{qt}{s}+p<1$, by Proposition \ref{P2.6}(i) and (\ref{410}) there exists $0<c<1$ such that
\[
cd(x)^s\leq u_i(x)\leq\frac{1}{c}d(x)^s
\]
in $\Omega$, $i=1,2$. Therefore  there exists a constant $C>1$ such that $Cu_1\geq u_2$ and $Cu_2\geq u_1$ in $\Omega$.

We claim that $u_1\geq u_2$ in $\Omega$. Supposing by contradiction, let
\[
\Gamma=\inf\{\gamma>1:\gamma u_1\geq u_2\text{ in }\Omega\}.
\]
By our assumption, we have $\Gamma>1$. From $\Gamma u_1\geq u_2$ in $\Omega$, it follows that
\[
(-\Delta)^{t} v_2=u_2^{-r}v_2^{-\theta}\geq\Gamma^{-r}u_1^{-r}v_2^{-\theta}
\]
in $\Omega$. Thus $v_1$ is a positive classical solution and $\Gamma^{\frac{r}{1+\theta}}v_2$ is a positive classical supersolution of
\[
\left\{
\begin{array}{rrll}
(-\Delta)^{t} w &=&u_1^{-r}w^{-\theta} & {\rm in} \ \ \Omega\\
w &=& 0  & {\rm in} \ \ \mathbb{R}^n\setminus\Omega
\end{array}
\right.,
\]
because $-r+\frac{r}{1+\theta}=-\theta\frac{r}{1+\theta}$. From the Proposition \ref{P2.1}, we obtain $v_1\leq\Gamma^{\frac{r}{1+\theta}} v_2$ in $\Omega$. Combining  the above estimate, we get
\[
(-\Delta)^{s} u_1=v_1^{-q}u_1^{-p}\geq\Gamma^{-\frac{qr}{1+\theta}}v_2^{-q}u_1^{-p}
\]
in $\Omega$. Therefore $u_2$ is a positive classical solution and $\Gamma^{\frac{qr}{(1+p)(1+\theta)}}u_1$ is a positive classical supersolution of
\[
\left\{
\begin{array}{rrll}
(-\Delta)^{s} z &=&v_2^{-q}z^{-p} & {\rm in} \ \ \Omega\\
z &=& 0  & {\rm in} \ \ \mathbb{R}^n\setminus\Omega
\end{array}
\right.,
\]
because $-\frac{qr}{1+\theta}+\frac{qr}{(1+p)(1+\theta)}=-p\frac{qr}{(1+p)(1+\theta)}$. By Proposition \ref{P2.1}, we conclude $u_2\leq\Gamma^{\frac{qr}{(1+p)(1+\theta)}} u_1$ in $\Omega$. Since $\Gamma>1$ and $\frac{qr}{(1+p)(1+\theta)}<1$, the above inequality contradicts the minimality of $\Gamma$. Then, $u_1\geq u_2$ in $\Omega$. Arguing similarly we conclude  $u_1\leq u_2$ in $\Omega$, so $u_1\equiv u_2$ which we obtain  $v_1\equiv v_2$. Thus, the system has a unique positive classical solution. This ends the proof of uniqueness.\fim \\


 \end{document}